\documentclass[12pt]{article}

\usepackage{amsfonts}
\usepackage{graphicx}
\usepackage[french]{babel}

\newtheorem{lemma}{Lemme}

\title{Le logo du CNRS est-il convexe ?}

\author{Didier Henrion$^{1,2}$}

\begin{document}

\maketitle

\footnotetext[1]{LAAS-CNRS, Universit\'e de Toulouse, France}
\footnotetext[2]{Facult\'e de G\'enie Electrique,
Universit\'e Technique Tch\`eque de Prague, R\'epublique Tch\`eque}

\addtocounter{footnote}{2}

\begin{abstract}
En octobre 2008 le CNRS adopte un nouveau logotype \`a la forme
arrondie. Nous \'etudions la repr\'esentation math\'ematique
de cette forme, et en particulier sa convexit\'e.
\end{abstract}

\section{Introduction}

\begin{figure}[h!]
\begin{center}
\includegraphics[width=\textwidth]{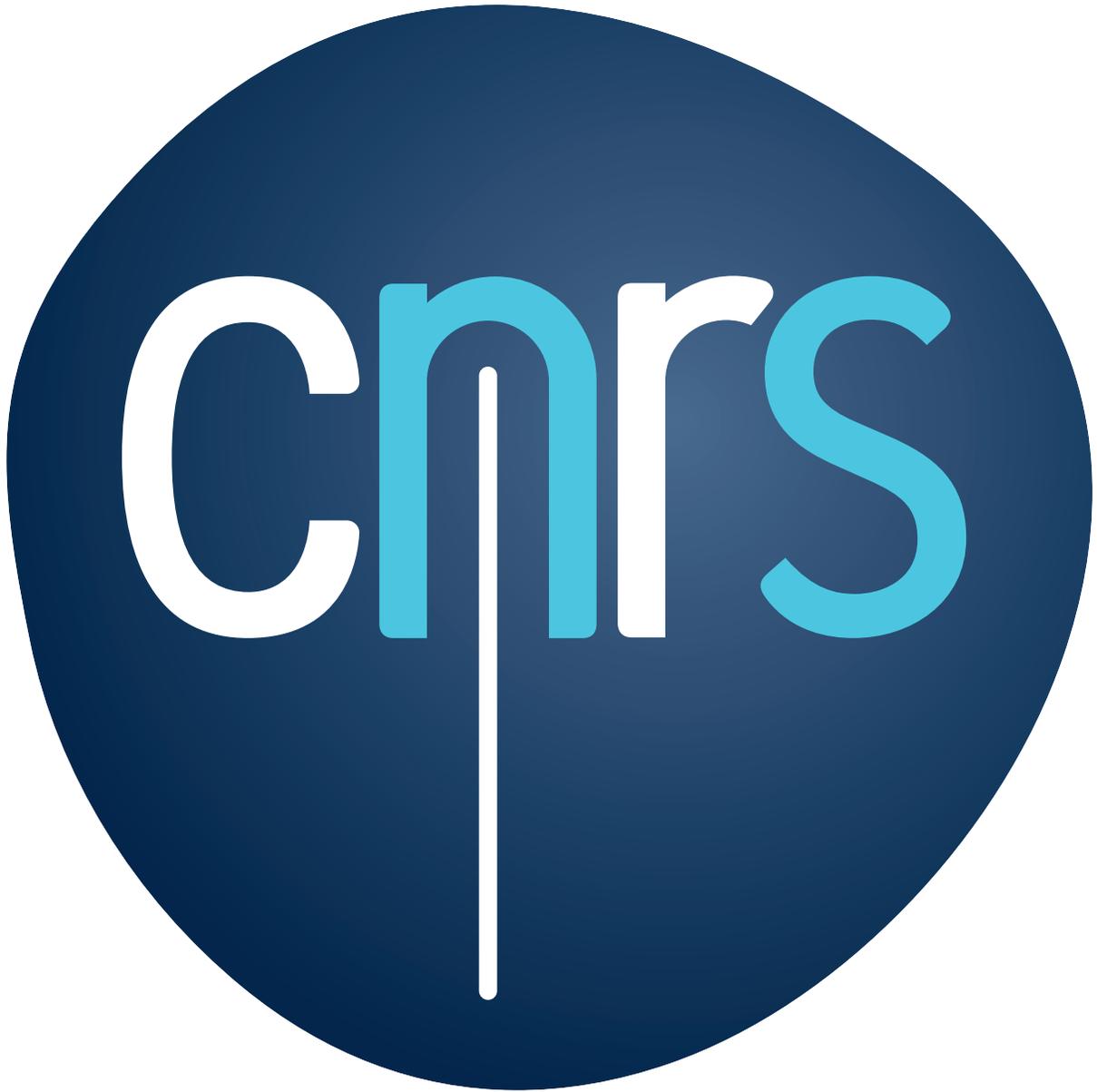}\\
\caption{Logotype du CNRS \`a la forme arrondie.
\label{logo}}
\end{center}
\end{figure}

En octobre 2008 le CNRS adopte une nouvelle identit\'e visuelle,
avec un logotype radicalement modifi\'e, caract\'eris\'e en ces
termes \cite{journal}: ``Les lettres de son sigle [..] sont
inscrites dans une forme arrondie et compacte qui exprime
l'unit\'e que conservera le nouveau CNRS. Pas tout \`a fait
ronde, novatrice, la forme du logo figure le processus m\^eme
de la recherche, toujours en devenir, et \'evoque la mati\`ere
mise \`a la disposition de nos chercheurs par notre plan\`ete.
Une mati\`ere mall\'eable, pr\^ete \`a se livrer aux expertises
de la recherche scientifique comme la motte de terre glaise
dans les mains du sculpteur'', voir la figure \ref{logo}.
Pour une interpr\'etation d'un autre registre, dans le contexte
politique actuel, voir \'egalement \cite{slr}.

\begin{figure}[h!]
\begin{tabular}{cc}
\includegraphics[width=.5\textwidth]{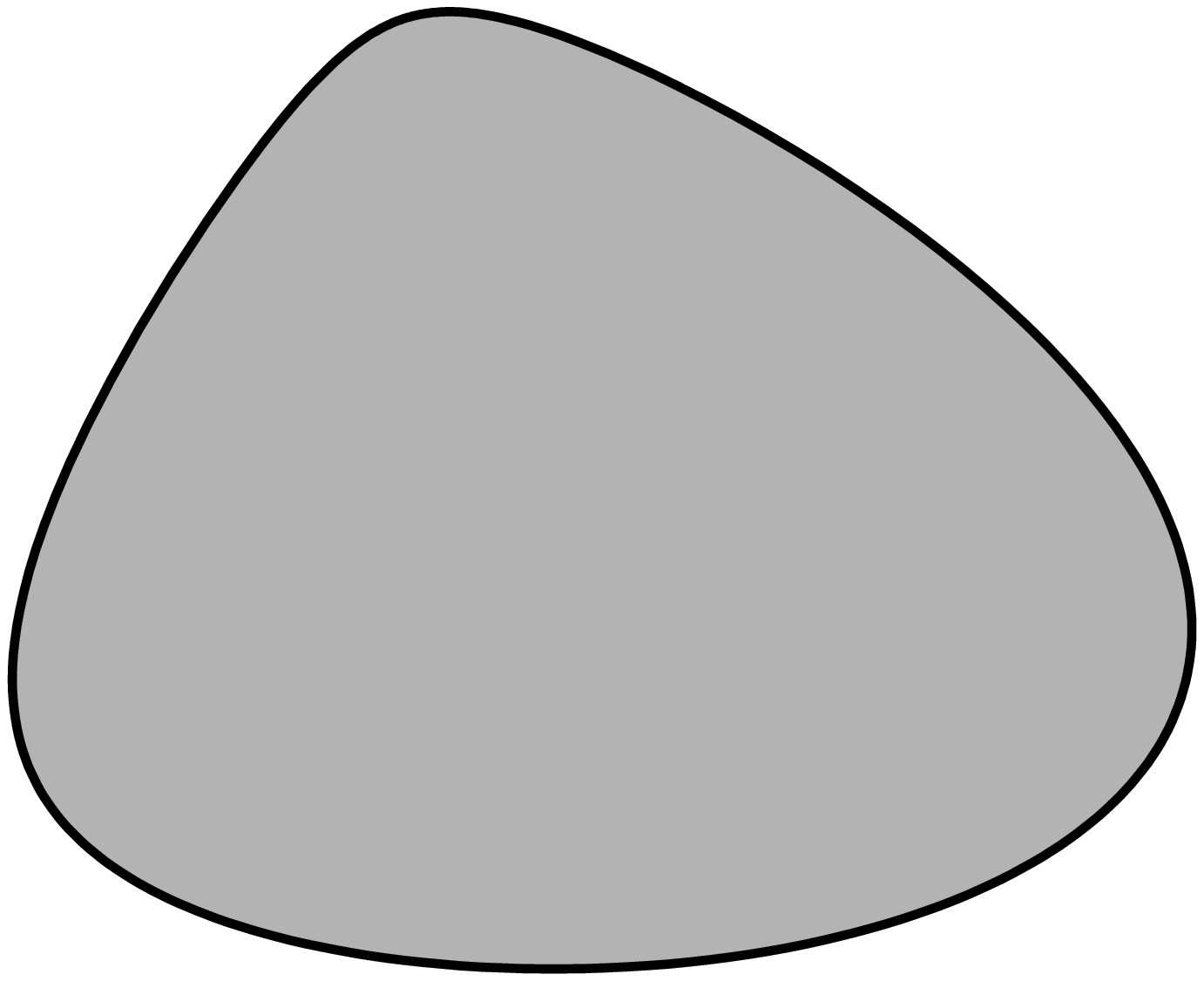}
&
\includegraphics[width=.5\textwidth]{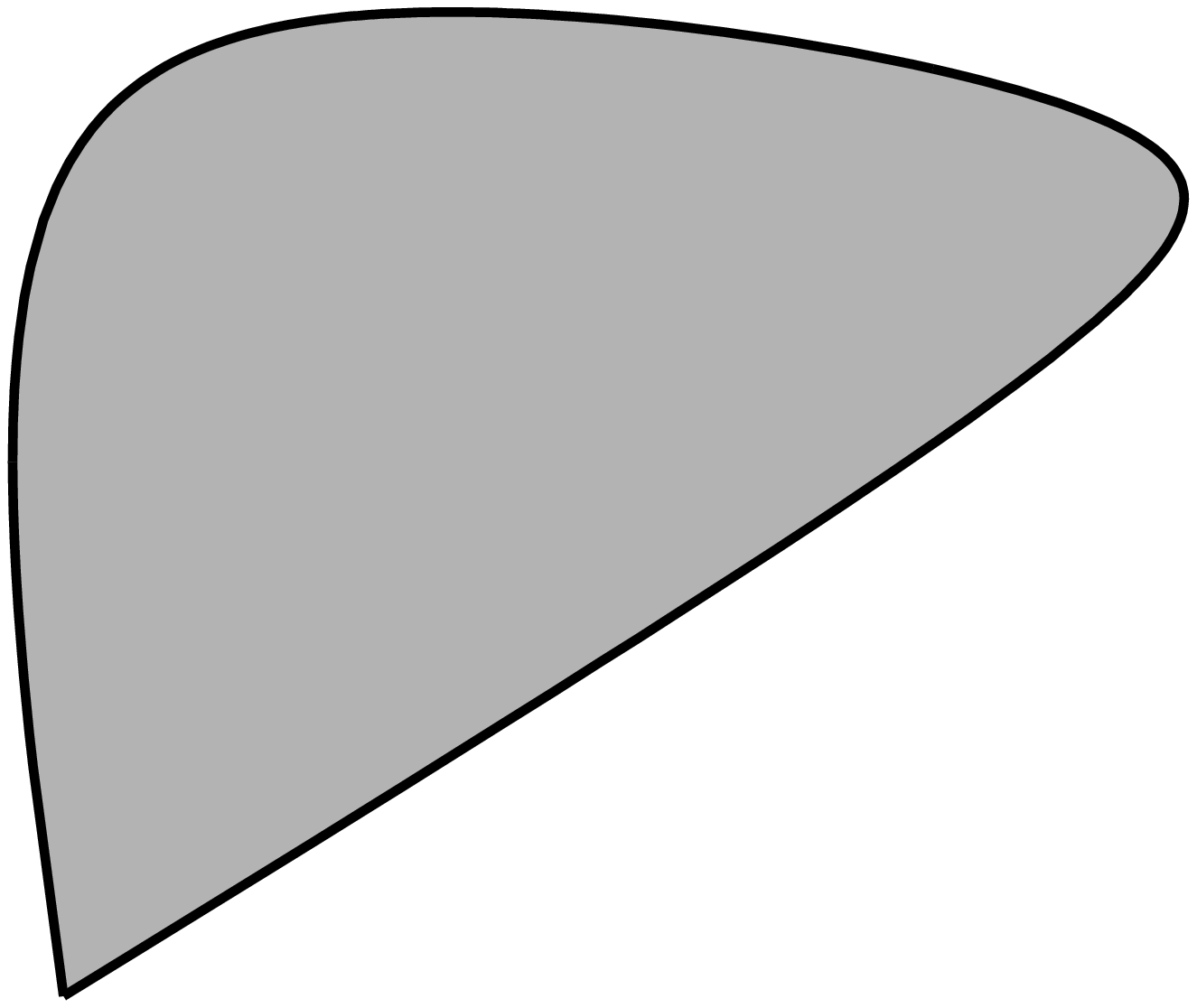}
\\
\includegraphics[width=.5\textwidth]{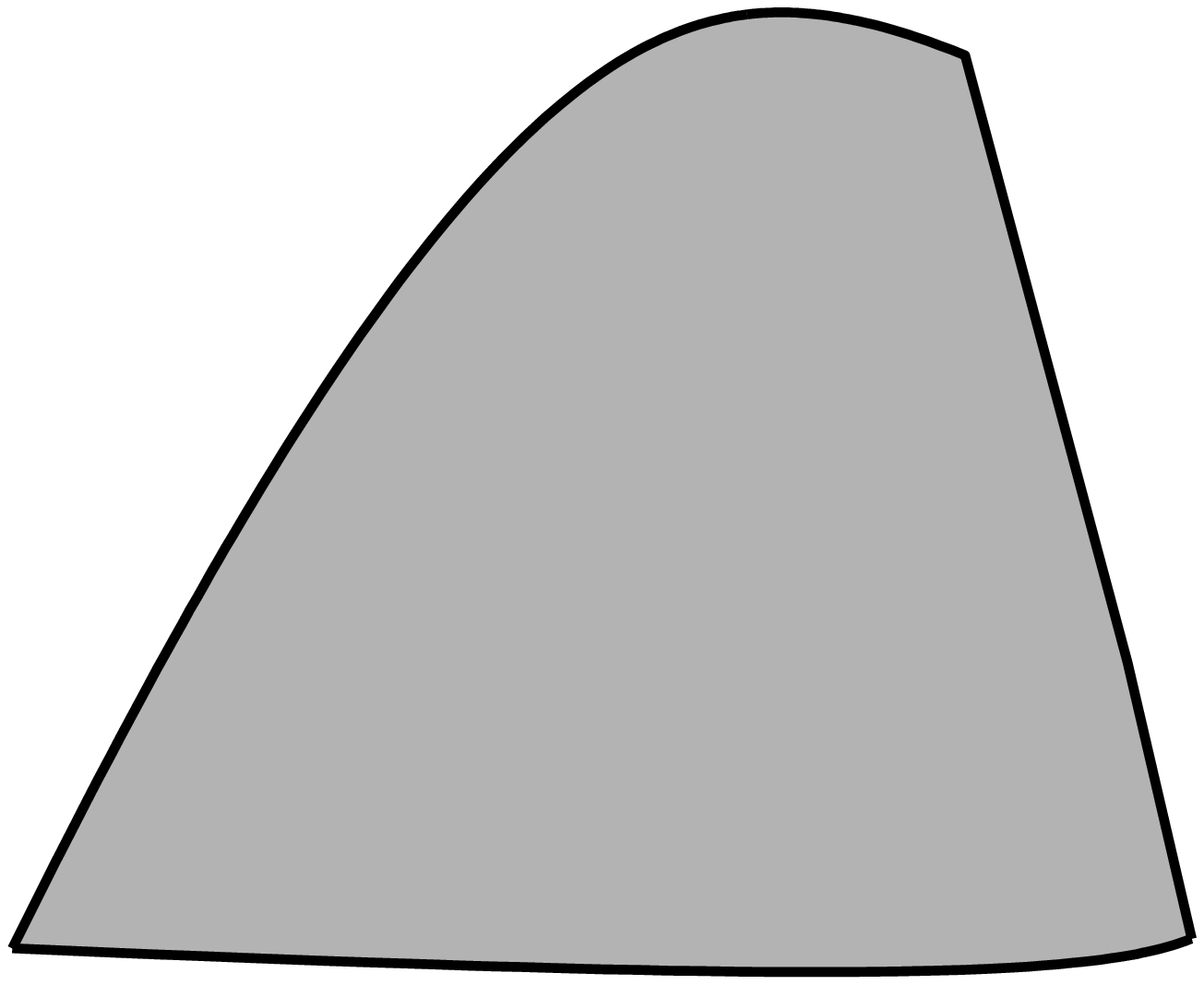}
&
\includegraphics[width=.5\textwidth]{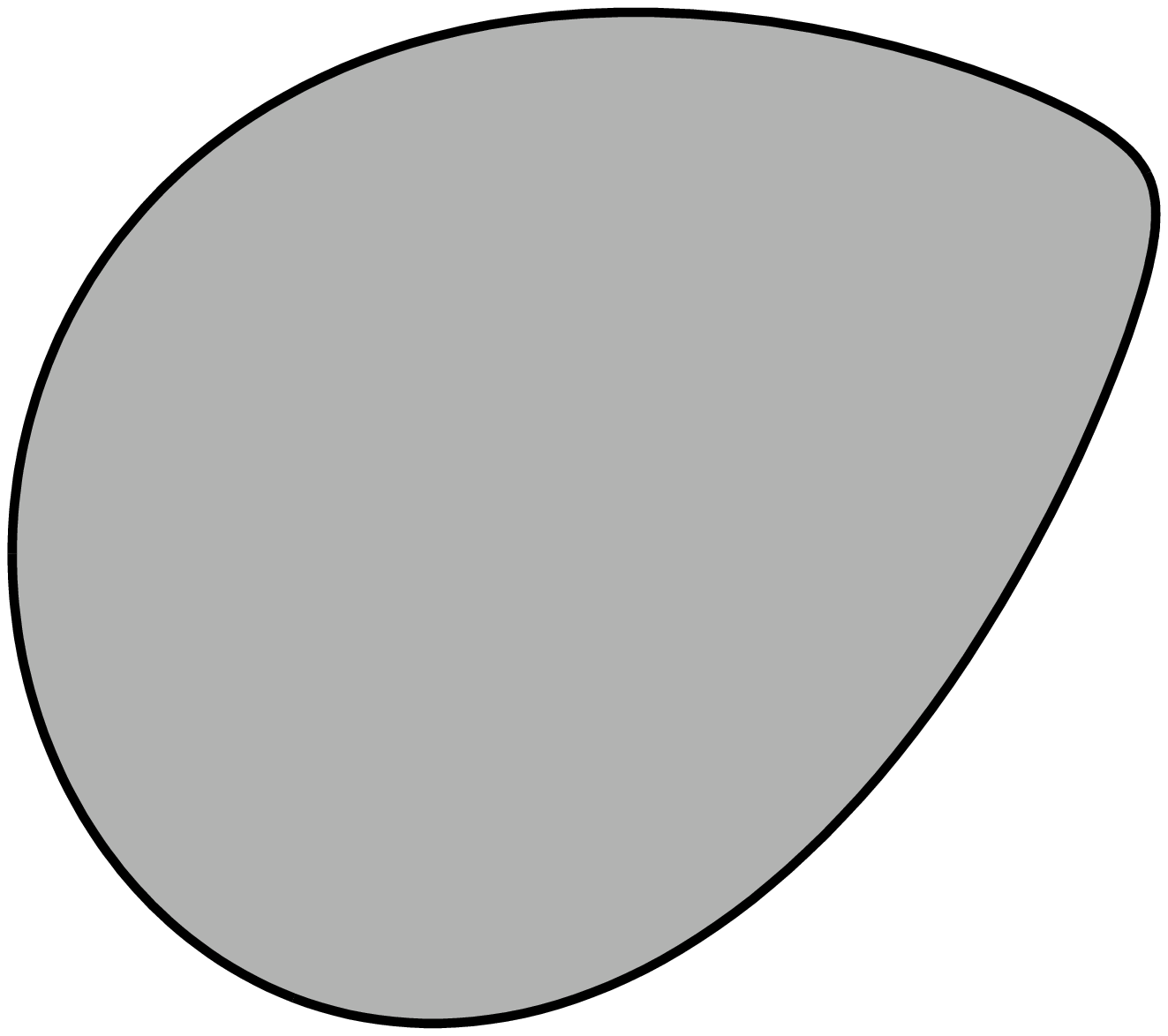}
\end{tabular}
\caption{Formes g\'eom\'etriques convexes.
\label{lmi}}
\end{figure}

En tant qu'utilisateur de techniques de programmation math\'ematique 
et d'optimisation \cite{jbhu}, j'ai \'et\'e souvent
confront\'e \`a des formes g\'eom\'etriques similaires, voir par exemple
la figure \ref{lmi}.
En optimisation, il est particuli\`erement souhaitable que ces
ensembles soient {\bf convexes}, voir par exemple \cite{bn}
ou \cite{boyd} pour plus de d\'etails.
Un ensemble est dit convexe lorsque,
pour tout couple de points choisis dans l'ensemble, le segment
reliant les points y est enti\`erement contenu. Toutes les
formes de la figure \ref{lmi} sont convexes. Un exemple
classique d'ensemble non-convexe est la lune ou croissant
apparaissant, entre autres, sur le drapeau de la Turquie.

La question qui m'est venue \`a l'esprit \`a la vision
de la figure \ref{logo} est donc la suivante :
\begin{center}
{\bf Le logo du CNRS est-il convexe ?}
\end{center}

\section{Courbes de B\'ezier}

La courbe du logotype a \'et\'e dessin\'ee par un graphiste
\`a l'aide du logiciel Adobe Illustrator, d\'evelopp\'e par
l'entreprise am\'ericaine Adobe Systems. Ce logiciel de
cr\'eation graphique permet de g\'en\'erer des images
vectorielles consistu\'ees de courbes \`a l'aspect
modulable gr\^ace \`a l'outil ``plume'' permettant
de placer des points d'ancrage et des tangentes.
Un des avantages des images vectorielles est qu'elles
ne sont pas d\'ependantes de la r\'esolution, c'est-\`a-dire
qu'elles ne perdent pas en qualit\'e si on les agrandit.

Dans Adobe Illustrator l'\'el\'ement de base
est la courbe de B\'ezier cubique, impl\'ement\'e via la
commande {\tt curveto} du language PostScript
\cite[page 564]{postscript}, le format natif des documents
g\'en\'er\'es par les logiciels d'Adobe. Les courbes de B\'ezier
sont des courbes polynomiales param\'etriques utilis\'ees
\`a l'origine au d\'ebut des ann\'ees 1960
par l'ing\'enieur fran\c{c}ais Pierre
B\'ezier pour la conception de pi\`eces d'automobile
chez Renault. Une repr\'esentation efficace de ces
courbes a \'et\'e propos\'ee par le math\'ematicien
et physicien fran\c{c}ais Paul de Faget de Casteljau,
alors qu'il travaillait pour Citro\"en.

\begin{figure}[h!]
\begin{tabular}{cc}
\includegraphics[width=.5\textwidth]{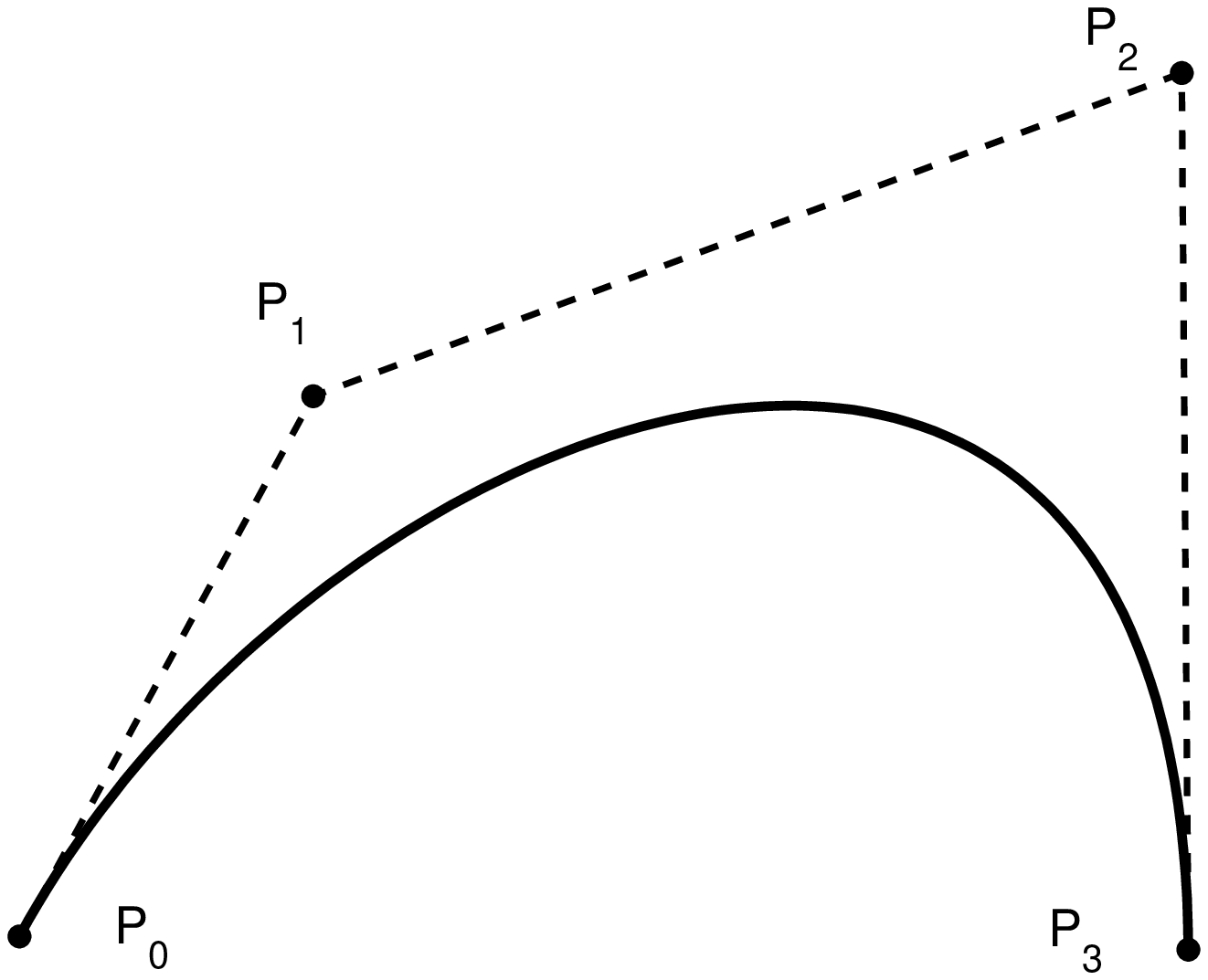}
&
\includegraphics[width=.5\textwidth]{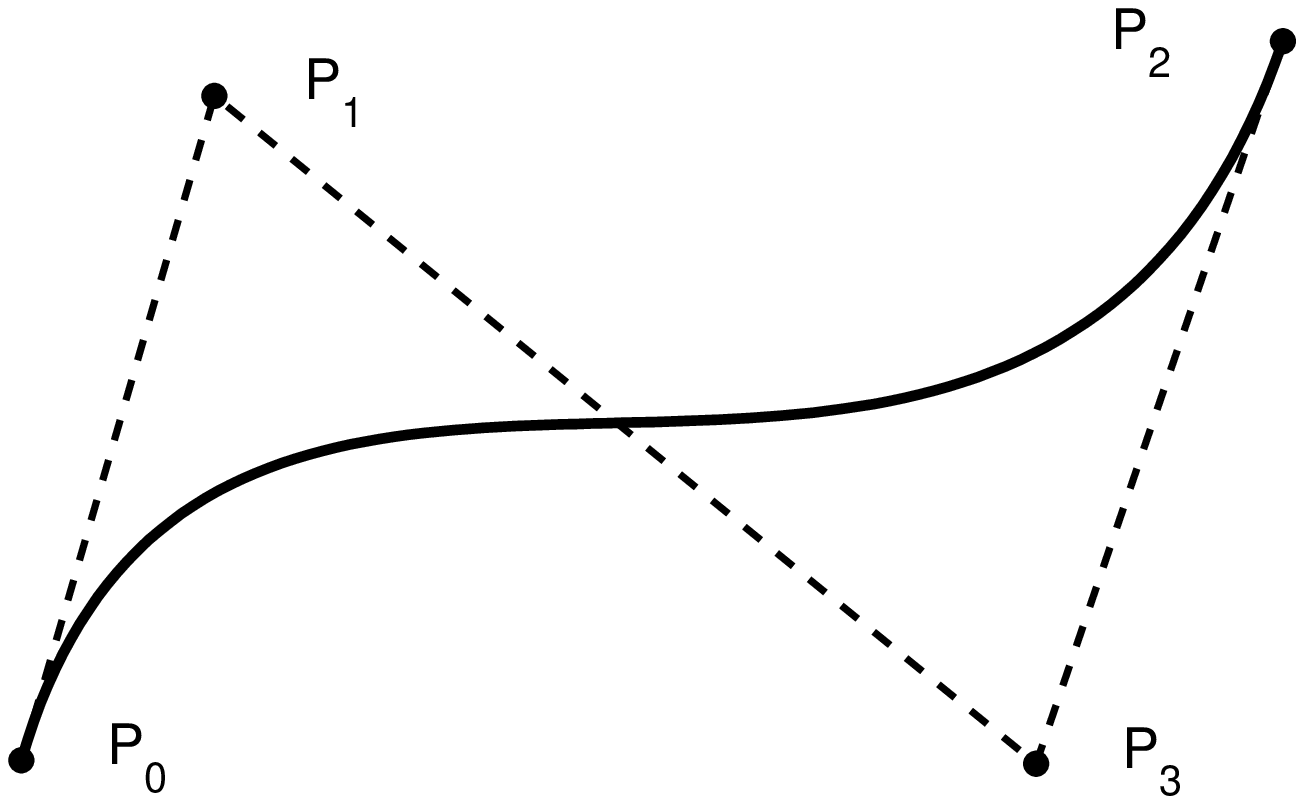}
\\
\includegraphics[width=.5\textwidth]{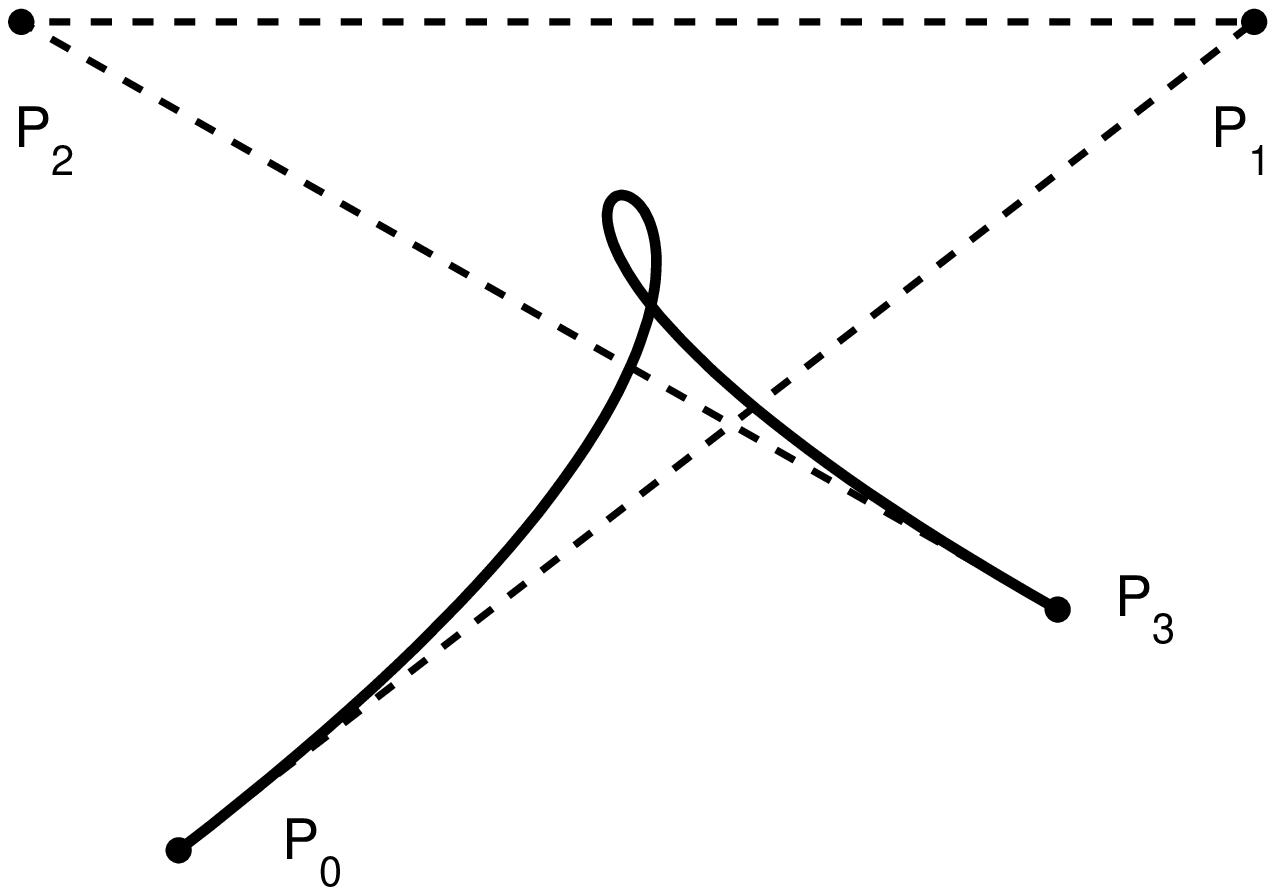}
&
\includegraphics[width=.5\textwidth]{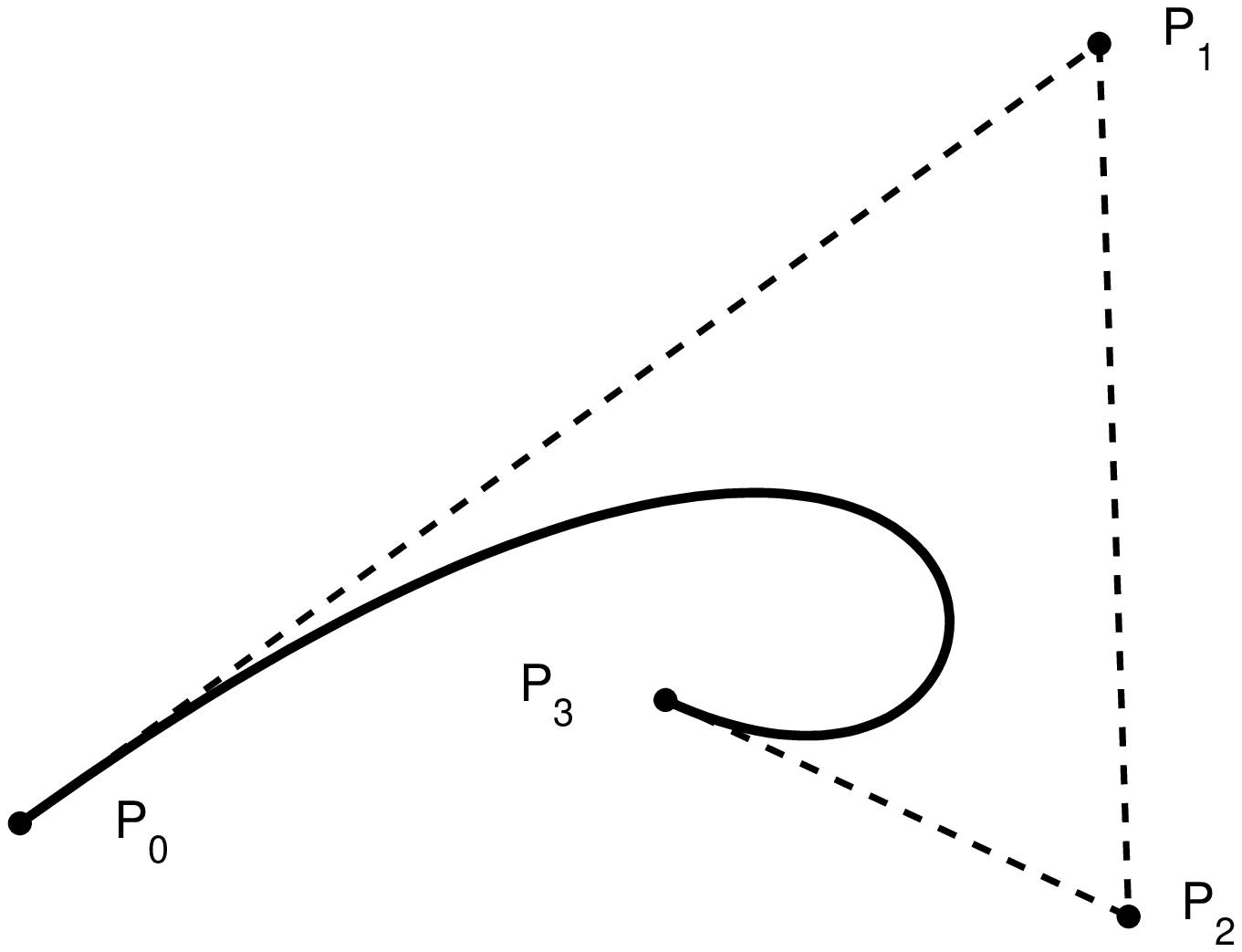}
\end{tabular}
\caption{Courbes de B\'ezier cubiques.
\label{bez}}
\end{figure}

Une courbe de B\'ezier cubique est d\'efinie par quatre
points de contr\^ole $P_0$, $P_1$, $P_2$ et $P_3$.
La courbe se trace en partant du point $P_0$ dans la
direction du point $P_1$, et en arrivant au point $P_3$
dans une direction venant du point $P_2$. En g\'en\'eral, la courbe
ne passe ni par $P_1$ ni par $P_2$, voir la figure
\ref{bez} pour quelques illustrations, et par exemple
\cite[Chapitre 4]{mortenson} pour un traitement rigoureux.

\begin{figure}[h!]
\begin{center}
\includegraphics[width=\textwidth]{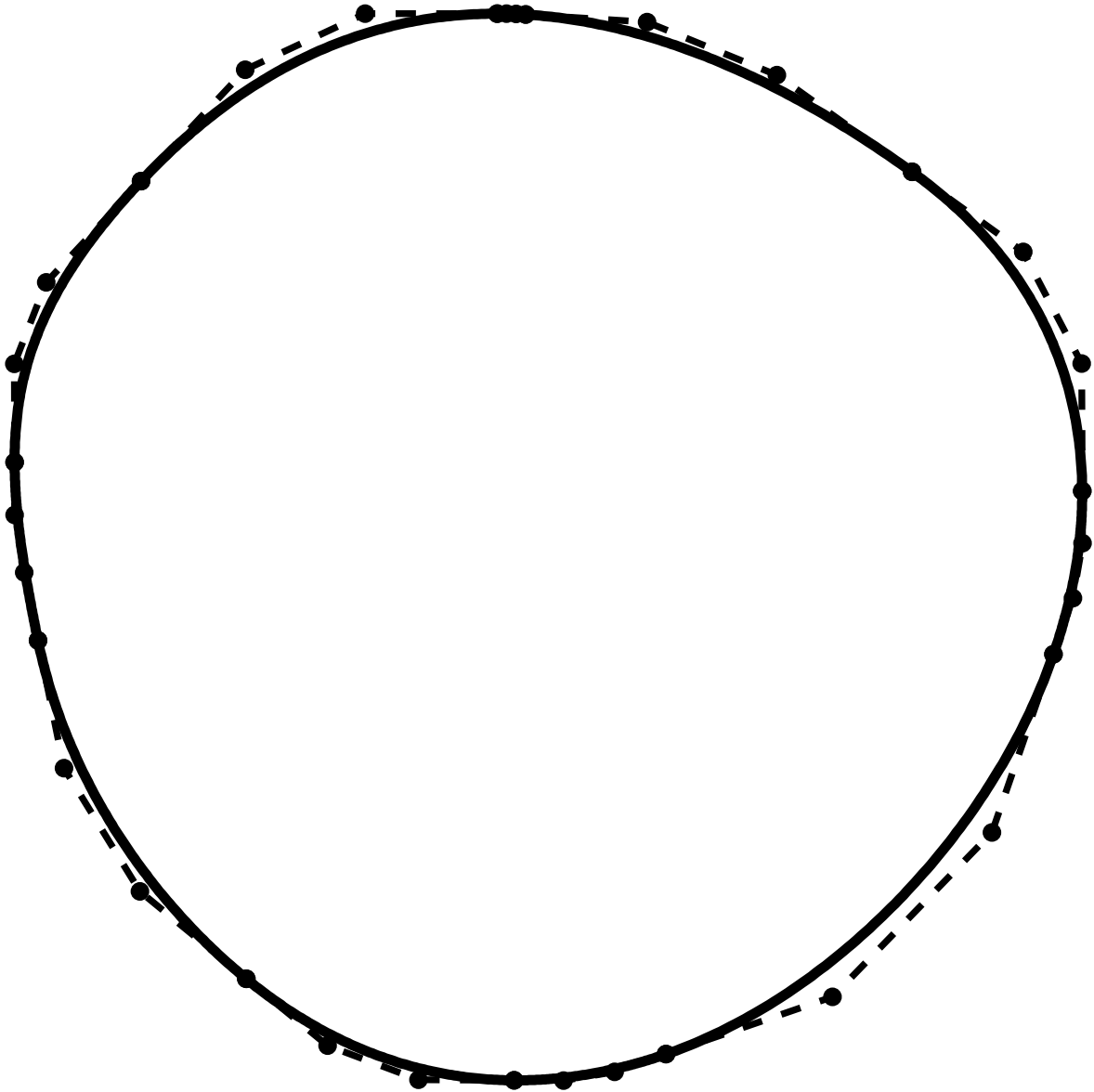}\\
\caption{Forme du logo constitu\'ee de 11 courbes de
B\'ezier cubiques et 33 points de contr\^ole.
\label{bezlogo}}
\end{center}
\end{figure}

L'exploration du fichier source PostScript du logotype du CNRS
montre que la forme est g\'en\'er\'ee par 11 instructions
{\tt curveto}, c'est-\`a-dire 11 courbes de B\'ezier cubiques
cons\'ecutives avec 33 points de contr\^ole, voir la
figure \ref{bezlogo}.

\section{Convexit\'e du logotype}

Afin d'\'etudier la convexit\'e du logotype du CNRS,
il faut donc \'etudier la convexit\'e des
courbes de B\'ezier qui le composent.

Etant donn\'es deux points $P_i$ et $P_j$, notons
$[P_i,\:P_j]$ le segment les reliant. Etant donn\'ee
une s\'equence ordonn\'ee de points $P_0,\:P_1,\:P_2,\ldots P_N$, notons
$[P_0,\:P_1,\:P_2,\ldots P_N,\:P_0]$ le polygone form\'e
par les segments $[P_0,\:P_1]$, $[P_1,\:P_2]$, $\ldots$,
et $[P_N,\:P_0]$.

Notons $C^0$ la courbe de B\'ezier
de points de contr\^ole $P_i$, $i=0,1,2,3$.
Notons $B^0$ l'ensemble d\'elimit\'e par $C^0$
et le segment $[P_0,\:P_3]$.
Notons $P^0$ le polygone $[P_0,\:P_1,\:P_2,\:P_3,\:P_0]$,
un quadrilat\`ere pouvant \^etre non-convexe, voir
la figure \ref{bez}.

\begin{lemma}
$B^0$ est convexe si et seulement si $P^0$ est convexe.
\end{lemma}

\begin{figure}[h!]
\begin{center}
\includegraphics[width=\textwidth]{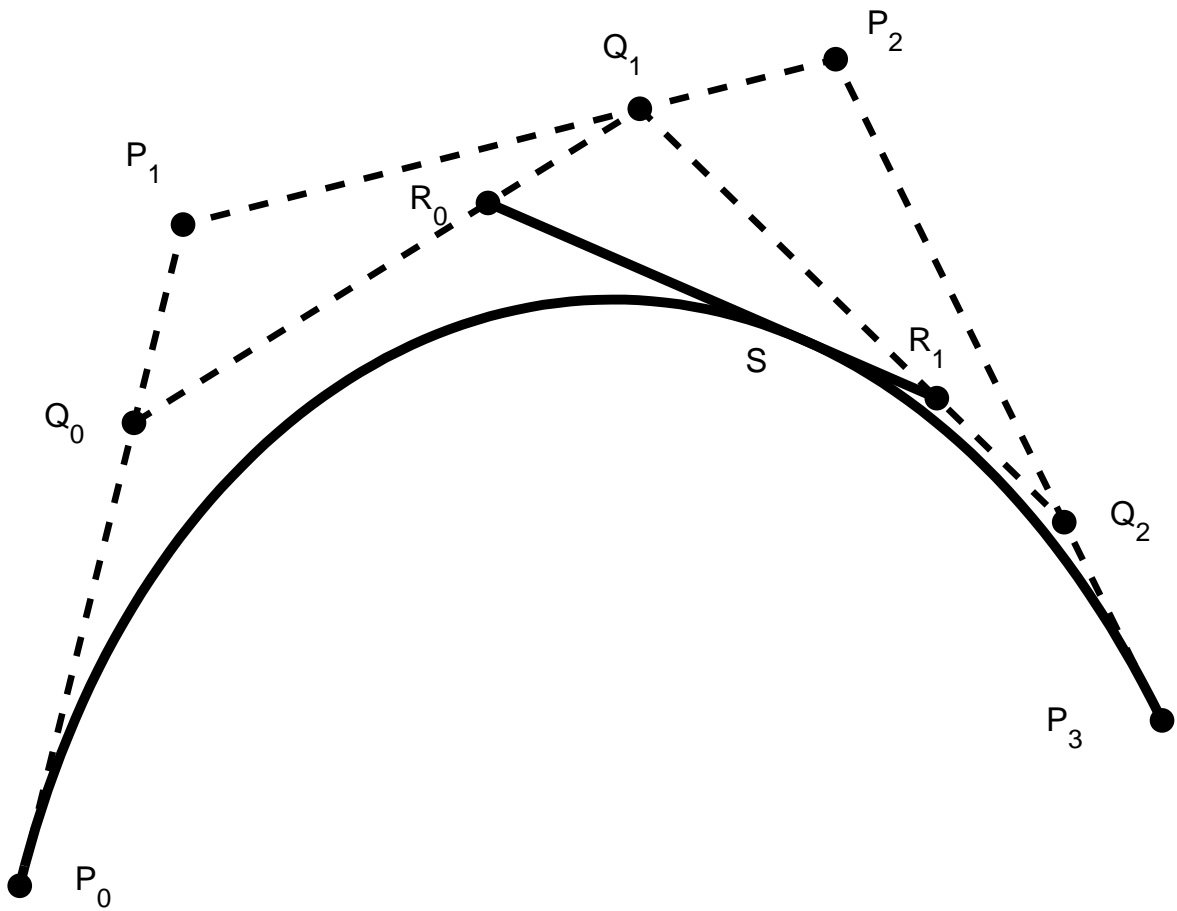}\\
\caption{Courbe de B\'ezier convexe et sa tangente
(traits pleins) construite par l'algorithme de de Casteljau.
\label{convexe}}
\end{center}
\end{figure}

{\bf Preuve:} 
Elle d\'ecoule de l'algorithme r\'ecursif
de construction propos\'e par de Casteljau et
d\'ecrit par exemple dans \cite[Section 4.1]{mortenson},
voir la figure \ref{convexe}. Pour un param\`etre $u \in [0,\:1]$
on d\'efinit les points $Q_i = uP_i+(1-u)P_{i+1}$, $i=0,1,2$
puis les points $R_i = uQ_i+(1-u)Q_{i+1}$, $i=0,1$.
Le segment $[R_0,\:R_1]$ est alors tangent \`a la courbe 
$C^0$ au point $S=uR_0+(1-u)R_1$. La courbe
est constitu\'ee de l'ensemble de ces points de tangence
lorsque $u$ varie entre $0$ et $1$. Si le polygone
$P^0$ est convexe, alors par construction le segment tangent
ne coupe jamais la courbe, quelque soit $u$, ce qui
implique la convexit\'e de l'ensemble $B^0$.
Invers\'ement si $B^0$ est convexe, alors le segment tangent
ne coupe jamais $B^0$, et donc le quadrilat\`ere $P^0$ est convexe.
$\Box$

Pour $k=1,\ldots,N$ notons \`a pr\'esent $C^k$ la courbe
de B\'ezier num\'ero $k$, de points de contr\^ole
$P^k_i$, $i=0,\ldots,3$. Supposons que les $N$ courbes
soient reli\'ees de telle sorte que
$P^k_3=P^{k+1}_0$ et $P^N_3=P^1_0$, et notons $B$
l'ensemble d\'elimit\'e par la courbe ferm\'ee
r\'esultante. Finalement, notons $P$ le polygone
$[P^1_0,\:P^1_1,\:P^1_2,\:P^1_3,\:P^2_1,\:P^2_2,\:P^2_3,\ldots,P^N_2,\:P^1_0]$.

\begin{lemma}
$B$ est convexe si et seulement si $P$ est convexe.
\end{lemma}

\begin{figure}[h!]
\begin{center}
\includegraphics[width=\textwidth]{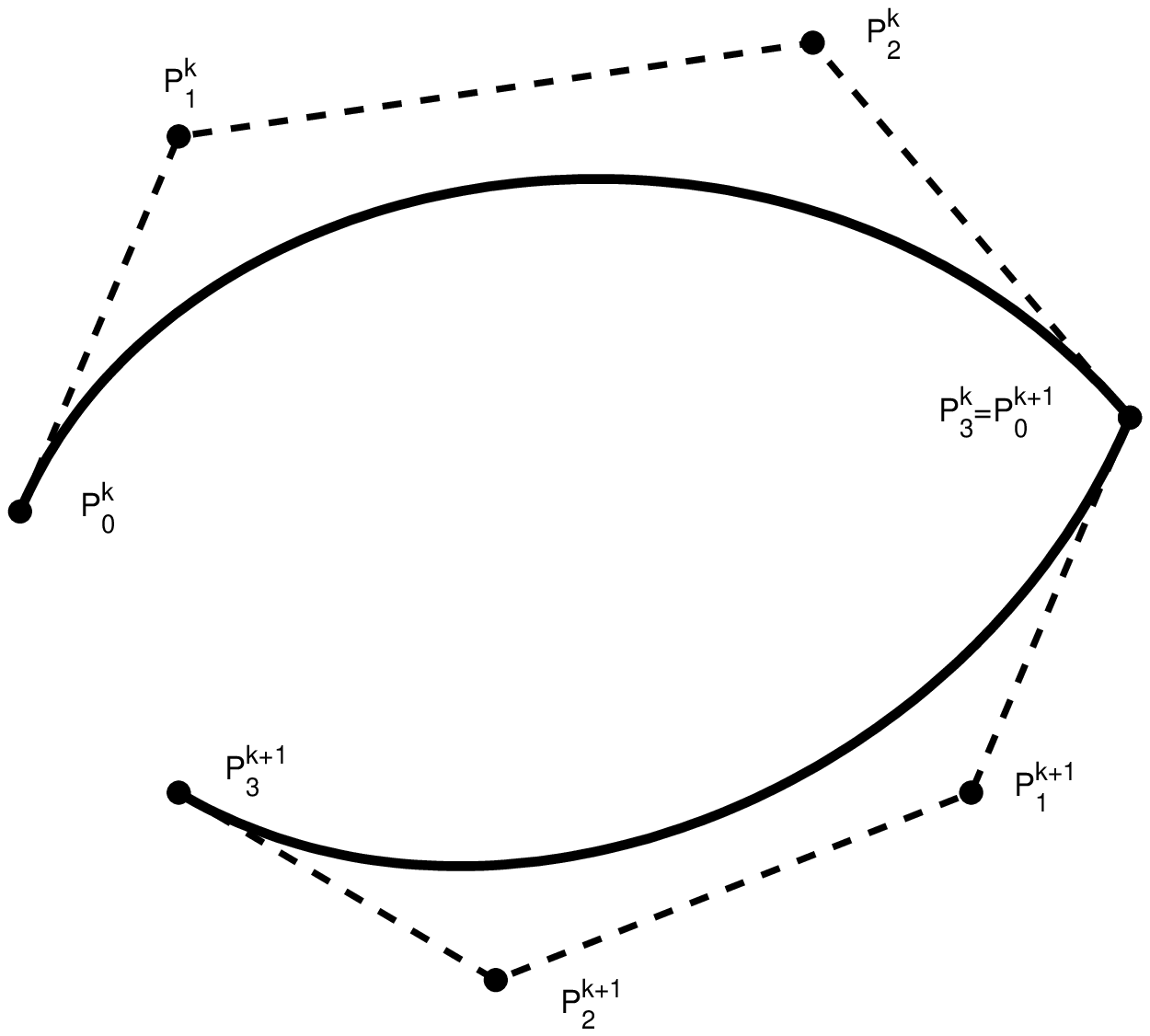}\\
\caption{Liaison convexe non-lisse.
\label{convexe1}}
\end{center}
\end{figure}

\begin{figure}[h!]
\begin{center}
\includegraphics[width=\textwidth]{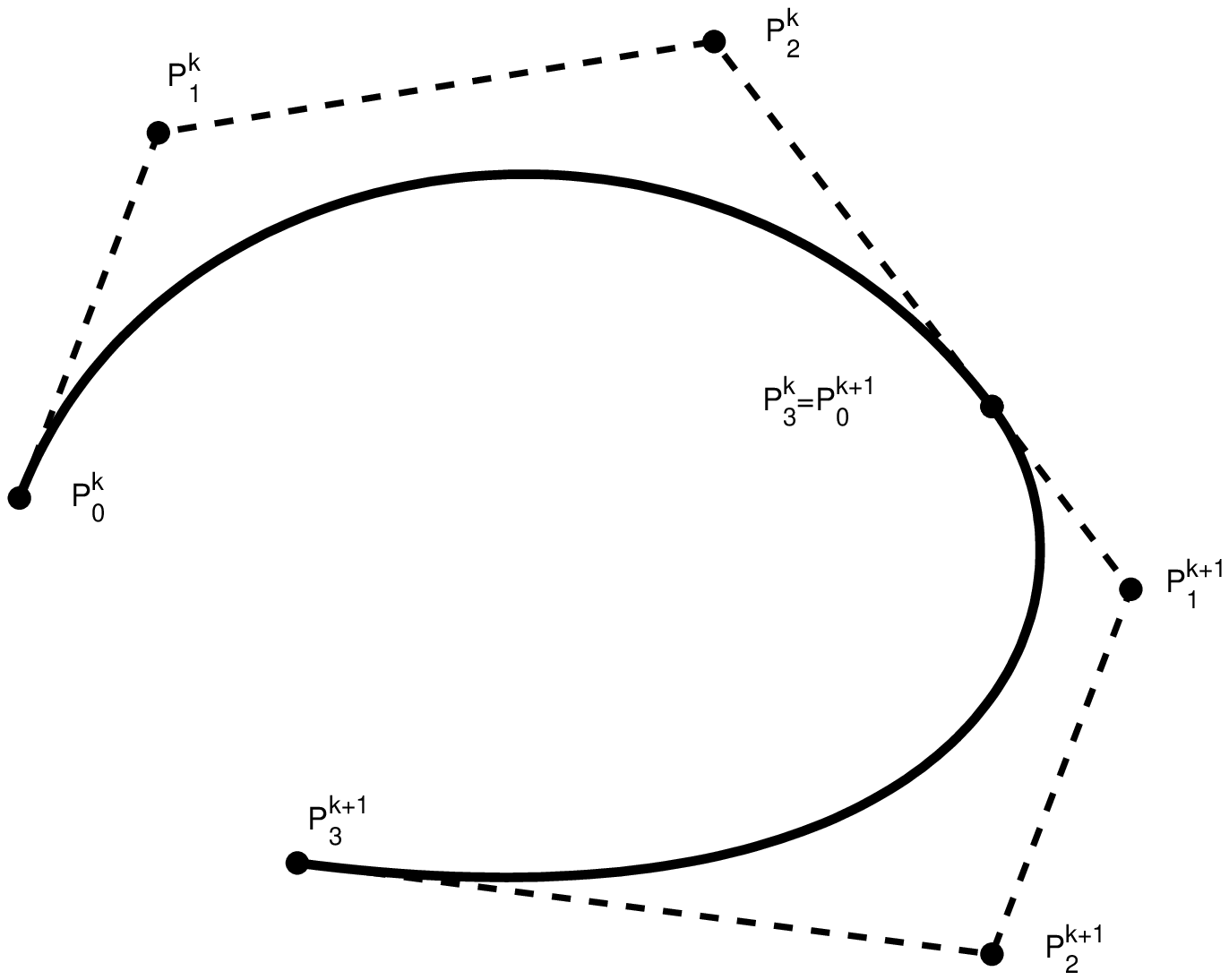}\\
\caption{Liaison convexe lisse.
\label{convexe2}}
\end{center}
\end{figure}

\begin{figure}[h!]
\begin{center}
\includegraphics[width=\textwidth]{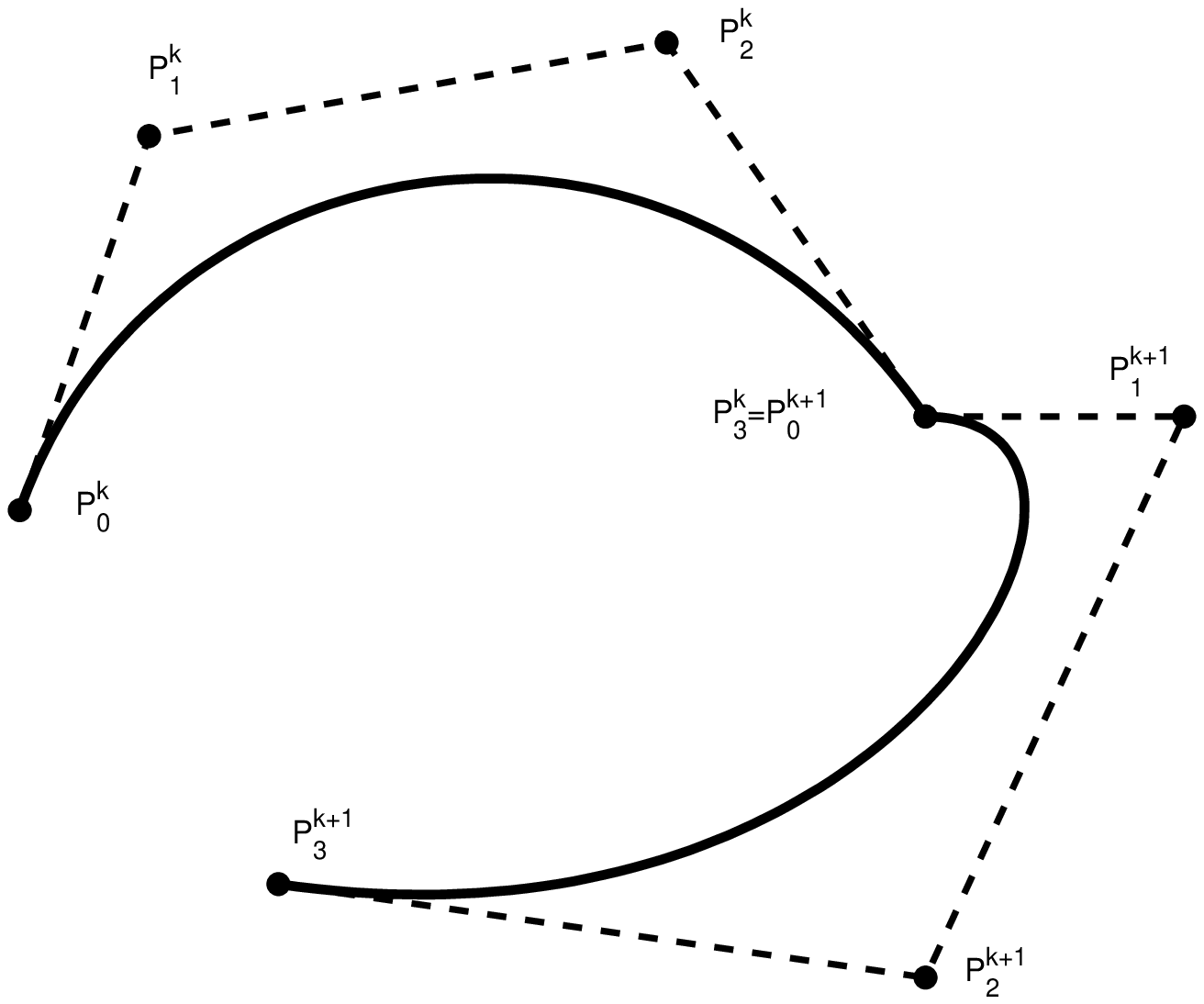}\\
\caption{Liaison non-convexe non-lisse.
\label{convexe3}}
\end{center}
\end{figure}

{\bf Preuve:}
Au point de liaison entre deux courbes convexes $C^k$ et $C^{k+1}$,
trois configurations sont possibles :
\begin{enumerate}
\item le segment $[P^k_2,\:P^{k+1}_1]$ coupe les courbes
$C^k$ et $C^{k+1}$, auquel cas l'ensemble d\'elimit\'e
par les courbes et le segment $[P^k_0,\:P^{k+1}_3]$
est convexe mais non-lisse au point de liaison,
voir la figure \ref{convexe1};
\item le segment est tangent aux deux courbes, auquel cas
les points $P^k_2$, $P^k_3=P^{k+1}_0$ et $P^{k+1}_1$
sont align\'es et la liaison entre les courbes est
lisse, voir la figure \ref{convexe2};
\item le segment ne coupe pas les courbes, auquel cas
l'ensemble d\'elimit\'e est non-convexe et non-lisse au
point de liaison, voir la figure \ref{convexe3}.
\end{enumerate}
La convexit\'e de $P$ implique la convexit\'e de
chaque quadrilat\`ere de contr\^ole, et
donc la convexit\'e de chaque cubique $C^k$. De plus,
la convexit\'e de $P$ implique que l'on ne
peut pas se trouver dans la troisi\`eme configuration
ci-dessus, et donc que l'ensemble $B$
d\'elimit\'e par les courbes est convexe.
Invers\'ement, la convexit\'e de $B$ implique
la convexit\'e de chaque courbe $C^k$, de leurs
liaisons et donc du polygone $P$.$\Box$

\begin{figure}[h!]
\begin{center}
\includegraphics[width=\textwidth]{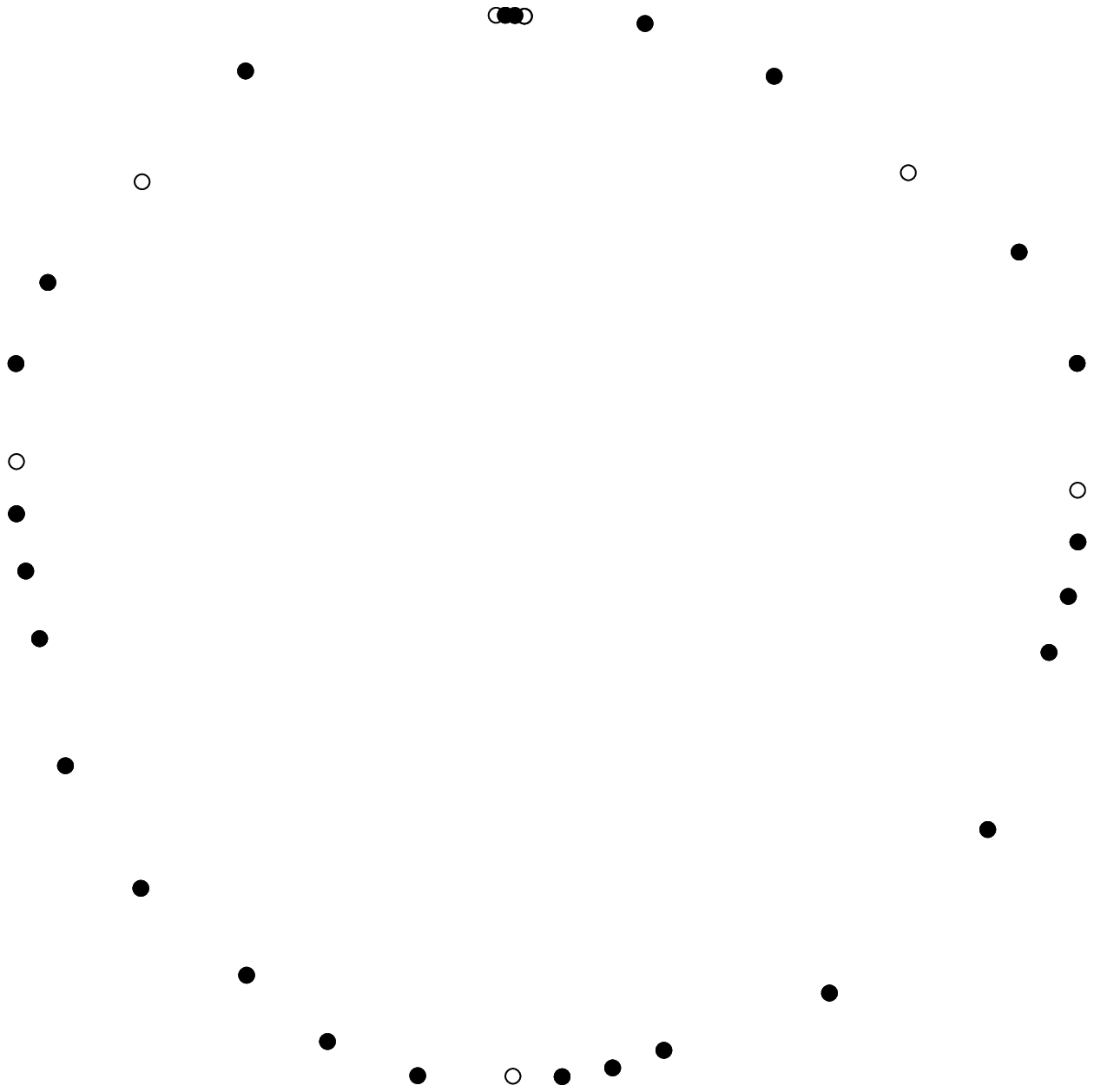}\\
\caption{Points de contr\^ole de l'enveloppe convexe
(disques noirs) et points de liaison responsables
de la non-convexit\'e (disques blancs).
\label{sommets}}
\end{center}
\end{figure}

Pour r\'epondre \`a la question de la convexit\'e du
logo du CNRS, il faut donc \'etudier la convexit\'e
de son polygone de contr\^ole. Ce polygone comporte
33 points. Pour savoir s'il est convexe, j'ai
d\'etermin\'e quels sont les sommets de son enveloppe
convexe (le plus petit ensemble convexe contenant
le polygone), \`a l'aide de l'algorithme d\'ecrit dans
\cite{qhull} et interfac\'e dans le logiciel de
calcul scientifique Matlab. Il s'av\`ere que seulement
26 points sur 33 g\'en\`erent l'enveloppe convexe,
voir la figure \ref{sommets}, d'o\`u la conclusion :
\begin{center}
{\bf Le logo du CNRS n'est pas convexe.}
\end{center}
Je m'abstiendrai de tirer une quelconque morale
de cette observation, ou d'essayer d'en donner
une interpr\'etation s\'emiologique comme dans \cite{slr}.

Les 11 quadrilat\`eres de contr\^ole des courbes de B\'ezier
sont tous convexes, mais 7 points de liaison
sont dispos\'es comme sur la figure \ref{convexe3}.
Cependant, en pratique, les points du logo sont
presque align\'es,
et l'effet de non-convexit\'e est beaucoup plus
faible que sur la figure \ref{convexe3}. Ainsi,
il ne semble pas possible de d\'etecter \`a l'oeil nu
la non-convexit\'e et la nature non-lisse de la courbe,
voir la figure \ref{bezlogo}.

L'expression math\'ematique exacte de la fronti\`ere
du logo du CNRS est celle d'une courbe cubique
convexe par morceaux, mais globalement non-convexe.
En d\'epla\c{c}ant l\'eg\`erement les 7 points
de contr\^ole il est possible d'assurer des
liaisons lisses et convexes, tout en pr\'eservant
l'aspect global du logo. Une autre possibilit\'e
consiste \`a g\'en\'erer une courbe de B\'ezier
convexe de degr\'e 25 en ne gardant que les 26 points
de contr\^ole de l'enveloppe convexe.

\section*{Remerciements}

Je remercie cordialement Anne-Solweig Gremillet
de la Direction de la Communication du CNRS pour
m'avoir procur\'e tr\`es rapidement le logotype
du CNRS sous format \'electronique, ainsi que
Dominique Daurat du Service de Documentation-Edition
du LAAS-CNRS pour son aide avec le format du fichier.

\begin{thebibliography}{XX}

\bibitem{postscript}
Adobe Systems Inc.
``PostScript language reference'',
3\`eme \'edition.
Addison-Wesley, 1999.

\bibitem{qhull}
C. B. Barber, D. P. Dobkin, H. T. Huhdanpaa
``The quickhull algorithm for convex hulls'',
ACM Trans. Math. Software, 22(4), pp. 469-483, 1996.

\bibitem{bn}
A. Ben-Tal, A. Nemirovski
``Lectures on modern convex optimization'',
SIAM, 2001.

\bibitem{boyd}
S. P. Boyd, L. Vandenberghe
``Convex optimization'',
Cambridge University Press,
2005.

\bibitem{jbhu}
J. B. Hiriart-Urruty ``L'optimisation'',
Collection ``Que sais-je ?'', Presses Universitaires de France,
1996.

\bibitem{mortenson}
M. E. Mortenson
``Geometric modeling'',
2\`eme \'edition.
Wiley, 1997.

\bibitem{slr}
E. Saint-James
``Le nouveau logo du CNRS, petite analyse s\'emiologique'',
post\'e sur le site de l'Association ``Sauvons la recherche''
le 7 octobre 2008.

\bibitem{journal}
C. Zeitoun
``Une nouvelle identit\'e visuelle pour le CNRS'',
Le Journal du CNRS, 225, p. 34, octobre 2008.
\end{thebibliography}
\end{document}